\documentclass{primus}

\usepackage{amsmath,amssymb}
\usepackage{color}
\usepackage{graphicx}
\usepackage{hyperref}

\title{Logistics of Mathematical Modeling-Focused Projects}
\author{R. Corban Harwood\\
Department of Mathematics\\
George Fox University\\
414 North Meridian Street, Newberg, Oregon, United States\\
rharwood@georgefox.edu\\
(503) 554-2737}

\newcommand{\AmSLaTeX}{$\cal A$\kern-.1667em\lower.5ex\hbox{$\cal
M$}\kern-.125em $\cal S$-\LaTeX}

\begin{document}


\makePtitlepage
\makePtitle

\begin{abstract}
This article addresses the logistics of implementing projects in an undergraduate mathematics class and is intended both for new instructors and for instructors who have had negative experiences implementing projects in the past. Project implementation is given for both lower and upper division mathematics courses with an emphasis on mathematical modeling and data collection. Projects provide tangible connections to course content which can motivate students to learn at a deeper level. Logistical pitfalls and insights are highlighted as well as descriptions of several key implementation resources. Effective assessment tools, which allowed me to smoothly adjust to student feedback, are demonstrated for a sample class.
As I smoothed the transition into each project and guided students through the use of the technology, their negative feedback on projects decreased and more students noted how the projects had enhanced their understanding of the course topics. Best practices learned over the years are given along with project summaries and sample topics. These projects were implemented at a small liberal arts university, but advice is given to extend them to larger classes for broader use.
\end{abstract}

\listkeywords

\section{INTRODUCTION}

One of my main joys in teaching is helping students draw connections between what they know and what they are learning. Projects can do just that. To effectively engage students in this way, I incorporate projects in most of my classes. Such projects are most naturally done in applied courses, but can be integrated into theoretical courses as well. 
My most effective projects have covered the full spectrum of mathematical modeling where students are involved in data collection, processing, development of the model equations, and evaluation of their model based upon the results and outside data.

\subsection{Learning Through Projects} \label{sec:literature}
A review of the literature on project-based learning and problem-focused group work has shown that these nontraditional teaching techniques support student learning in a variety of ways. Project-based learning is a teaching methodology that utilizes student-centered projects to facilitate learning, while problem-focused group work utilizes motivational problems to set the stage for the group's investigation. Implementing projects increases student engagement in class, collaboration outside of class, and gives students new perspectives of the material. Long-term studies by Hestenes and Epstein have shown that interactive teaching styles result in significantly higher understanding of concepts \cite{Hestenes,Epstein}. In classes with roughly 40 students, Mergendoller showed that short-term projects improved student attitude and achievement \cite{Mergendoller}. In larger classes with roughly 75 students, Olson showed that long-term problem-focused approaches similarly improved student learning and achievement \cite{Olson}. In developing class projects, core principles of project-based learning can be helpful guides, such as making group collaboration effective and motivating investigation with open-ended questions. Specifically, group collaboration that involves student choice, communication, writing, revision, and presentation is most effective at increasing student learning \cite{Linhart}. Additionally, group projects provide stimulating discussions and can spur ideas for individual research projects later on.



\subsection{The Modeling Perspective}
Modeling-focused projects are helpful in providing a big-picture perspective. Starting a course with such a project can motivate students immediately. For example, on the first day of my Mathematical Modeling course, I asked students to guess the number of hazelnuts which I could pack into the jar I held in my hand. Then, given only the empty jar, one hazelnut, a ruler, and ten minutes, the groups quickly estimated the number by comparing the volumes of a hazelnut and the jar. After a class discussion reflecting on the modeling process, students were prepared to research packing efficiency and come back to the next class with better estimates. The room buzzed with anticipation the following class day as I counted out hazelnuts and packed them into the jar. Some students groaned as they realized how far off their estimates were, while others sat with bated breath. After tightening the lid, shaking the jar and squeezing as many hazelnuts in as I could, the winning student shouted with delight as I ended exactly on her estimate. I wrapped up the project by discussing how much the modeling process improved the estimates, both in terms of a closer average and smaller spread. Our final discussion gave a great segue into statistical analysis.

This project exemplifies another core principle of project-based learning: students developing mastery and becoming self-directed learners \cite{Ambrose}. Since parameter estimation from collected data was incorporated into this project, students achieved a more realistic view of the world which helped them utilize their learning later on \cite{Winkel}. Mooney and Swift frame the modeling process as creating an idealized `Model World' through simplifying assumptions which are helpful and necessary due to available resources \cite{Mooney}. The solution to the idealized problem is then evaluated in the real world and the model is improved and solved again as necessary \cite{Bliss}. Starting with a quick visual guess, then volume measurements, and finally accounting for packing efficiency, my students could see the improvement with each successive resolution of their model. This simple activity demonstrates how and why we solve problems in the idealized Model World as well as how many assumptions surround even simple calculations. 

The modeling process can be reinforced by checking in with students during or after a project using questions such as {\it ``Why was this a good way to set up the real question mathematically?;" ``Is this a good mathematical technique for solving the problem?;''  ``How well can you trust your answer?;''} or {\it``What do you do with the solution?'' .} Not knowing how to solve the real problem exactly, students felt doubtful in giving estimated solutions to problems they framed themselves, but encouraging further justification, such as a logical argument developing their model and statistical tests of the data, helped them gain confidence. Facing such open-ended questions in a supportive learning environment helps prepare students for collaborative work in their jobs where they will have to explain and defend their results \cite{Winkel}.

\subsection{Outline}
The rest of this article is organized as follows. Implementation difficulties, utilized resources, and a generalization of best practices are shared in section \ref{sec:logistics}. Section \ref{sec:implementation} summarizes specific project implementations, chosen topics, and student feedback for classes in which I have implemented major projects. Section \ref{sec:appendix} discusses the main points of this article and reflects upon their impact on students.
Additionally, detailed project prompts, sizes of classes and groups utilized, and project scoring guides are included in the Appendix. Template files guiding these projects can be accessed online through this journal.
 
\section{LOGISTICS}\label{sec:logistics}
There are many logistical issues that can make implementing projects difficult for both professor and student. These can be minimized, as outlined in section \ref{sec:resources}, through a vareity of best practices discussed in section \ref{sec:bestpractices}.

\subsection{Logistical Issues}
The main logistical issue for implementing projects is the additional professor workload. Projects are an important way to engage students in a different learning style and encourage positive collaboration. Yet, projects can end up displacing important content and requiring much more time to prepare and grade. Further, it feels risky to dedicate an entire class period to group activities.

Also, issues with initializing group work need to be dealt with before the projects can be beneficial to student learning. Each student needs to know the content of the whole project, not just their individual component for the project to be an adequate learning tool. The grade should reflect the quality of the project as a whole and the individual's contribution. Unequal division of work is common when group members do not consciously divide up the work at the beginning. Also, groups who do not plan out their work end up completing the project at the last minute. 

\subsection{Resources}\label{sec:resources}
I have utilized several different tools in my classes, both in preparing projects and implementing them, to minimize the impact of these logistical issues. Databases of prepared projects and ideas save prep time. Collaboration tools and software lessen the difficulty in student coordination of the project. Course management tools house project work in one place for ease of group collaboration, in-class presentation, online submission, and grading. 

There are many resources available with prepared class projects. Toews gives a broad overview of how modeling can be used across the mathematics curriculum \cite{Toews}, while others cater to specific subjects such as calculus \cite{Karaali, MOSAIC}, differential equations \cite{Claus,SIMIODE}, numerical methods \cite{Cline}, and mathematical modeling with an emphasis on writing \cite{Linhart}. For its depth, I highlight {\it SIMIODE}, a teaching repository of modeling activites specifically for differential equations which is maintained by a community of educators \cite{SIMIODE}. These activities are peer-reviewed to provide clear instructions for student investigation and instructor facilitation. A similar community-supported repository specifically for calculus is {\it Project MOSAIC} \cite{MOSAIC}. Spurred by the 2013 Mathematics for Planet Earth initiative, the Center for Discrete Mathematics and Theoretical Computer Science at Rutgers University sponsored the development of several sustainability modules in the same vein but applied more broadly to calculus, differential equations, discrete mathematics, statistics, and liberal arts mathematics courses \cite{DIMACS}. Further, individual projects covering a wide spectrum of topics and classes can be found in the literature through journals such as PRIMUS.

Coordinated file management can aid in project implementation. A central database for web submission, display, and grading of projects has been a major time saver for me, even in a small class. File submission, forum setup, and wiki creation are supported by most course management systems (my university uses Moodle). A private wiki, accessed through a course management system like Moodle (\href{http://www.moodle.com}{www.moodle.com}) or separately like PBWorks (\href{http://www.pbworks.com}{www.pbworks.com}), is a web environment, similar to Wikipedia, where students can create and link together multiple web pages in a more controlled environment. Having students complete their project on the wiki keeps everything in one place to make it easier for me to grade. I can also visually check on group work progress and give directed reminders to those trailing behind. I often require students to self-report their individual contributions on the wiki to help them stay accountable and allow me to distribute points according to contribution efficiently. In addition, the fact that students are self-reporting their contributions encourages them to consciously divide up the work evenly at the beginning. The wiki environment also allows for quick  transition between oral presentations of multiple groups since they all link back to the page with the project's prompt. 

To help coordinate group work, I originally encouraged students to use mass emails to myself and their group, but now exclusively use a private wiki (through Moodle). Mass emailing is helpful because the emails can be sorted together and a record of the conversation is kept, but it does clutter my Inbox and email does not transfer mathematical work well. I now encourage students to share their work on the private wiki which allows posting of figures and mathematical typesetting in HTML. In addition to the benefits mentioned above, the wiki environment allows students to immediately see group updates in the project as a whole, add comments for revisions, and track changes to hold students accountable for contributing their share of the work. 

For collaborative data analysis, Google Sheets was effective in providing access to data collection, processing, and visualization. This saved me and the students the hassle of transferring files and made grading easier by keeping all of my students' work in the same document. Because every group worked on the same document, I created a seperate tab for each group with the first tab for instructions and example formatting. Both the wiki and Google Sheets are handy in group presentations as they decrease the amount of down time transferring files. Cloud storage systems such as Google Drive (drive.google.com) and Drop Box (dropbox.com) also provide a limited platform for file submission and group collaboration.

\subsection{Best Practices}\label{sec:bestpractices}
The number one advice for implementing new projects (or anything new) in your classes is to start small. If the technology needed for the project is new to the students, it is especially important to work out beforehand the amount of class demonstration time that is needed to train students. 

As a personal example, I first implemented one project (listed in section \ref{sec:PredPrey}) in just one class (Differential Equations with Linear Algebra). I chose this class and project for the ease of integration with the class content and the simplicity of working with the graphical interface of the Java applet (PPLANE) I had found for this project. Through discussions with students and course evaluations, I assessed the impact of the project and, as a result of positive feedback, built upon it's effectiveness with more projects.
 
What follows is a summary of best practices I have learned in designing effective assessments through student feedback, minimizing class preparation time, smoothly implementing the project, and extending the project to connect with areas outside of class. Note that the following comments are neither sufficiently exhaustive in their advice nor necessary to be implemented altogether.

\subsubsection{Assessment}
To know whether a change, such as a new project, in a class is beneficial, it is important to plan effective assessment tools for evaluating it. Student feedback in course evaluations have been somewhat helpful for me, especially when used repeatedly to compare iterations of a course, but getting timely and constructive feedback was often difficult. What follows is a detailed list of the assessment tools and refinements I used in a particular class to get clearer and more useful feedback from students about projects.

I first
created an anonymous midterm evaluation to guide students' reflection of their learning in terms of the course goals and learning components (such as projects). This feedback was much more focused than course evaluations, but sifting through all the comments was overwhelming. So, I categorized the feedback to track the overall impact of projects. 
Project-related feedback were grouped  by hindrance themes (H1-H3) and benefits (B1-B3). An ``X" in a semester marked when more than one comment related to this theme occurred. 
\begin{itemize}
\item[H1] Timing of projects with respect to lecture topics
\item[H2] Coordinating schedules for group work outside of class
\item[H3] Technical problems
\item[B1] Abstract topics became more tangible
\item[B2] Topics made more sense at a deeper level 
\item[B3] Good rhythm of individual homework and group projects 
\end{itemize}

\begin{table}

\begin{tabular}{lccccccc}
 & Fa12 & Sp13 & Fa13 & Sp14 & Fa14 & Sp15 & Fa15\\
Hindrances&&&&&&&\\
H1 &X&&&&&&\\
H2 &X&X&X&X&&&\\
H3 &&X&X&&X&&\\
&&&&&&&\\
Benefits&&&&&\\
B1 &X&&X&X&X&X&X\\
B2 &&&X&&X&X&X\\
B3 &&&&X&X&X&X\\
\end{tabular}
\caption{Assessment Tool Summarizing Repeated Occurrences of Project-Specific Themes in Student Feedback for Differential Equations with Linear Algebra from 2012-2015 (Sp=Spring, Fa=Fall)}\label{table:feedback} 
\end{table}



To demonstrate the utility of this assessment tool, the tenor of project-related feedback in sections of differential equations with linear algebra from Fall 2012 to Fall 2015 are summarized in Table \ref{table:feedback}. 
Tracking student feedback by semester, I noticed a positive correlation in the volume of benefit-themed comments vs. time and a negative correlation in the volume of hindrance-themed comments vs. time. As I repeated this class and had more experience implementing projects, the effect of logistical issues diminished and their benefits increased. This showed that the benefits on student learning increased once logistical issues were reduced. In sum, this assessment tool demonstrated that my implementation of projects in this class improved from Fall 2012 to Fall 2015.

Though comment summaries such as Table \ref{table:feedback} were helpful visualizations of my progress, I desired more quantitative information. Recently, I rewrote these open-ended questions on my midterm survey to have students share feedback by rating their level of agreement with clear statements about the impact of projects on their learning. I then supplemented my university's course evaluation with the following statements, to which student could respond: `Strongly Agree,' `Somewhat Agree,' `Neutral,' `Somewhat Disagree,' or `Strongly Disagree.'
\begin{enumerate}
\item Projects for this course added new perspectives on the content covered in homework.
\item Projects for this course helped me understand and engage with course concepts more fully

\end{enumerate}
In Spring 2016 course evaluations for differential equations with linear algrabra, students overwhelmingly agreed (89\% and 74\% in sections A and B) that these projects added new perspectives with similar agreement 
(84\% and 68\% in sections A and B) that these projects helped them understand and engage with the course concepts.

In sum, having clear and direct statements through which students can evaluate your course helps you make adjustments to successfully implement projects. 
Feedback specific to each sample project is listed in the project summaries in section \ref{sec:implementation}.

\subsubsection{Preparation}
To minimize class preparation time it is important to start with just one project that you have tried yourself, preferably prepared and tested by another source such as those mentioned in section \ref{sec:resources}. 

In preparing classes, I sought a good rhythm between teaching the needed skills and placing projects. I found a very structured project fit better at the beginning and a more open-ended project fit better in the latter part of the course. Embedding extensions of lecture content in the project helped with course cohesion as the homework problems fed into the project. I also found that projects two weeks or longer gave a thorough team-building experience through coordinated communication and scheduling of group work.

Incremental notification and implementation of projects has been key to merging project work with standard classwork. To keep from suddenly jarring students from the lecture mode, I found notifying them a week beforehand most helpful and had them form groups the class period before we started the project so they were ready to dive in together. 

When I assigned a project, I went over the expectations and demonstrated the technology students would be using whether it was a calculating web applet, wiki environment, or program they must compile themselves. Leading by example helped avoid many technical issues for students.

\subsubsection{Implementation}
To optimally implement a project, it is important to facilitate group work effectively.

I included opportunites for students to buy into a project through selecting their partners, their group role, and sometimes even their topic. To get students into groups they liked without the groups becoming too lopsided, I tried several different things. I first assigned groups myself (balanced in terms of performance in class so far), but found that students did not always work well together. To include student choice while making groups more diverse, I allowed students to pair up themselves and assigned groups from these pairs. Increasing student choice has led to groups bonding more quickly and working together more smoothly. Diversifying the groups in terms of major, gender, ethnicity, work ethic, etc. helped students in large classes get to know more of their classmates and break out of isolating cliques. For reduced grading, I personally have had most success building up to 12 groups per class in powers of two formed from these self-selected pairs: two for fewer than 24 students, four for 24-48 students, and eight for 48-96 students. At a small liberal arts university, my classes have ranged from 7 to 43 with a mean of 34, but at a public university earlier on in my career I found groups of eight (loosely paired groups of four) worked well for classes of 75 (Pre-Calculus) and 97 (Linear Algebra). Note that the projects in these large classes had a guided computer activity completed by the group outside of class instead of using an in-class activity to initiate the project.

It was important to remind students of their share of the workload and keep them accountable. I regularly checked in with groups to see if each student was aware of how his or her component fit into the project as a whole. Afterwards I had students describe his or her role in the group. Such evaluation of roles helped me monitor individual contributions. At times, I have asked students to evaluate their other group members but found this soured the project experience and I obtained better information by having group roles assigned at the beginning. My limited experience in assigning group roles was positive and it demonstrated better communication within groups. As one possible distribution in a group of four, I assigned roles of ``Leader" (distributes workload), ``Recorder" (records results of all components), ``Auditor" (checks work done for simple errors), and ``Spokesperson" (reports summary of project work).

\subsubsection{Connections Outside of Class}
To support student learning in the course as a whole, projects should enhance learning breadth and depth as well as the accessibility of content. Through the project requirements, I encouraged students to connect to other disciplines and to reference other researchers' work and data. To add depth, students should be accountable for their own work and communicate it well in written and oral forms. I used technology (such as a Moodle wiki and Google Sheets) to improve accessibility for students to collaborate, present visualizations of their work, and post reports online.

\section{PROJECTS}\label{sec:implementation}
Project summaries and student feedback are given in this section for three lower division and three upper division mathematics courses: liberal arts mathematics, discrete mathematics, differential equations, numerical methods, mathematical modeling, and advanced linear algebra. For open-ended projects, I kept a list of topics from which students chose for themselves. See section \ref{sec:ProjApp} for detailed project prompts and scoring guides. Ideas for these projects came from textbook resources and  the online community repositories, SIMIODE \cite{SIMIODE} and Project MOSAIC \cite{MOSAIC}, mentioned in section \ref{sec:resources}. When I needed to broaden a project from these sites, I often generalized parameters of the model, had students collect or analyze the dataset to be used, or extended the analysis of the solution.

\subsection{Liberal Arts Mathematics}
The liberal arts mathematics course at my university, {\it The World of Mathematics}, emphasizes applications of mathematical concepts in areas such as consumer finance, probability, and statistics. The {\tt Monthly Budget Spreadsheet.xlsx} template file can be accessed online and uploaded to Google Sheets as a collaborative class spreadsheet.

\textbf{Monthly Budget}
\begin{itemize}
\item Provided Resources
\begin{itemize}
\item Collaborative discussion forum (Moodle)
\item {\tt Monthly Budget Spreadsheet.xlsx} (Upload to Google Sheets)
\item Salary statistics (\href{http://www.salary.com}{www.salary.com})
\item Cost of living (\href{http://costofliving.salary.com}{costofliving.salary.com}) 
\item Mortgage payments (\href{http://www.zillow.com}{www.zillow.com})
\item Budget estimates (\href{http://www.learnvest.com}{www.learnvest.com})
\item TEDx talk on budgeting
(\href{http://www.tedxwallstreet.com/alexa-von-tobel-one-life-changing-class-you-never-took-2/#sthash.pfIh6s4v.dpbs}{www.tedxwallstreet.com/alexa-von-tobel-one-life-changing-class-you-never-took-2/})
\end{itemize}
\item Expectations
\begin{itemize}
\item Collect data for a chosen job, location, and house
\item Complete a budget and compare to the ideal 50-20-30 rule
\item Compare future value of retirement contributions for several initial amounts (simple sensitivity analysis)
\item Choose a best budget from these results and compare within group
\end{itemize}
\end{itemize}
Students were surprised at how short their ``paid vacation'' was in retirement, and overall liked how this project connected course topics to their lives.

\subsection{Discrete Mathematics}
Covering graph theory, analysis of algorithms, and various other topics, discrete mathematics brims with applications but is scattered over many topics. Thus, I chose one culminating graph theory project on game analysis spanning four weeks. In groups of four, students chose games from a list I provided and developed individual strategies which they tested against each other. To space this long-term project, I had groups teach the class their game in week two, before completing their analysis and presenting to the class in week four of the four-week project.

\textbf{Game Analysis}
\begin{itemize}
\item Resources
\begin{itemize}
\item Collaborative wiki (Moodle)
\end{itemize}
\item Expectations
\begin{itemize}
\item Choose a board/card game to summarize and analyze
\item Develop and test individual strategies against each other
\item Develop a board state evaluation function and demonstrate it on a limited game tree
\item Complete wiki report and ten minute presentation
\end{itemize}
\end{itemize}

This project was a highlight of the course for most students and I noticed an improvement in exam scores afterwards compared to before the project.

\subsection{Differential Equations with Linear Algebra}
Because differential equations are one of the most natural ways to express a dynamical model, this class provided ample opportunity for modeling projects. I have assigned either three or four 2-week projects in this 15-week course. The three projects I have consistently used are shown below. The {\tt Common Cold Outbreak Spreadsheet.xlsx} template file can be accessed online  and uploaded to Google Sheets as a collaborative class spreadsheet.

\subsubsection{Common Cold Outbreak Model}
\begin{itemize}
\item Resources
\begin{itemize}
\item Dorm floor plan; bag of beans
\item {\tt Common Cold Outbreak.xlsx} (Upload to Google Sheets)
\item Collaborative private wiki (Moodle)
\end{itemize}
\item Expectations
\begin{itemize}
\item Simulate spread of common cold by shaking beans onto floor plan and tracking infected rooms
\item Develop model and use data to estimate parameters
\item Complete wiki report with five minute presentation
\end{itemize}
\end{itemize}

Students were enthusiastic about this hands-on activity and appreciated seeing how a differential equation model developed.

\subsubsection{Forming a Nutritional Meal Plan}
\begin{itemize}
\item Resources
\begin{itemize}
\item Nutritional information on food from store or online
\item {\tt Matlab Nutrition Plan.m} (accessed online)
\item Collaborative private wiki (Moodle)
\end{itemize}
\item Expectations
\begin{itemize}
\item Collect nutritional information on foods to form a balanced day of meals
\item Form an initial meal plan, evaluate the overconsumption error, and compare it to the naive matrix solution and improved linear programming solution.
\item Report your findings in the wiki and present them to the class
\end{itemize}
\end{itemize}

Students found this application of a large matrix very insightful, especially since they constructed it from data they had collected. They enjoyed the creative aspect in forming their meal plan but were surprised at how hard it was to balance the fifteen nutritional areas we were tracking. Seeing matrix methods beyond the traditional Gaussian Elimination was eye-opening.

\subsubsection{Predator-Prey Model}\label{sec:PredPrey}
\begin{itemize}
\item Resources
\begin{itemize}
\item Collaborative private wiki (Moodle)
\item PPLANE applet (\href{http://math.rice.edu/~dfield/dfpp.html}{math.rice.edu/$\sim$dfield/dfpp.html})
\end{itemize}
\item Expectations
\begin{itemize}
\item Research qualitative and quantitative data on realistic predator and prey species
\item Develop predator-prey model and estimate parameters
\item Analyze model theoretically and via PPLANE simulation
\item Submit wiki report with ten minute presentation
\end{itemize}
\end{itemize}

Most groups chose standard predator-prey pairs such as orca-sea lion or cheetah-gazelle and enjoyed sharing biological/ecological tidbits learned in studying such animal behavior from around the globe. Every semester, however, I have a couple groups with more imaginative selections such as zombies-humans, dragons-sheep, and Sith-Jedi. These groups collected data from novels or games on which they based their parameters and compared their results. While less quantitative in nature, such a project provided a creative outlet for students. After this project in Fall 2013, a student who had been dragging her feet through the class so far exclaimed, ``I never knew math could be fun!'' It has been exciting to see student attitudes impacted so positively by projects.

\subsection{Numerical Methods}
Following several short technique-learning projects in numerical methods, the final five-week project gave students the flexibility of applying methods covered in class to topics of their choice. 
A schedule of weekly update reports have kept groups accountable along with a few group meetings with me to keep them on-track.

\textbf{Final Project}
\begin{itemize}
\item Resources
	\begin{itemize}
	\item Mathematical software (Matlab)
	\item Mathematical typesetting software (\href{http://www.lyx.org}{www.lyx.org})
	\item Collaborative forum (Moodle)
	\item {LyXReport.lyx} and {LyXPresentation.lyx} templates (accessed online)
	\end{itemize}
\item Expectations
	\begin{itemize}
	\item Develop short-term and long-term goals in numerically solving and/or analyzing a chosen problem
	\item Evaluate and update goals each week
	\item ten-page report typeset in TeX
	\item 15 minute presentation with one-page handout
	\end{itemize}
\end{itemize} 

Students appreciated the flexibility of this project. Several enhanced projects from other courses, such as {\it Simulating a Helical Solenoid} and {\it Identifying Predictors of Voting Registration}, while others extended topics introduced in this class such {\it Computing Land Area from a Sample of Longitude and Latitude} and {\it Investigating Fractal Behavior of Various Iteration Functions}. 

\subsection{Mathematical Modeling}
The entire mathematical modeling course was set up as a series of two-week group projects. 
Below is a summary of two of these projects extended from the textbook \cite{Mooney}: estimating hazelnuts in a jar and modeling the spread of the common cold on campus using an SIR model. The {\tt Common Cold SIR Spreadsheet.xlsx} can be accessed online  and uploaded to Google Sheets as a collaborative class spreadsheet.

\subsubsection{Hazelnut Estimation}
\begin{itemize}
\item Resources 
	\begin{itemize}
	\item Ruler, hazelnut, and jar dimensions 
	\end{itemize}
\item Expectations
	\begin{itemize}
	\item Measure volume of hazelnut and estimate how many fit in the given jar
	\item Research optimal packing efficiencies and improve initial group estimate
	\item Submit report summarizing chosen simplifying assumptions, solution method, and final estimate
	\end{itemize}
\end{itemize} 

Students enjoyed this project as a way to see how many individual decisions and various approaches can go into solving such a simple problem.

\subsubsection{Common Cold SIR Model}
\begin{itemize}
\item Resources
	\begin{itemize}
	\item Mathematical Software (Matlab)
	\item {\tt Common Cold SIR Spreadsheet.xlsx} (Upload to Google Sheets) 
	\end{itemize}
\item Expectations
	\begin{itemize}
	\item Use an anonymous survey of university students to track signs of a cold every day for a week (40-50 per group)
	\item Record data in the spreadsheet and compute average transmission probabilities
	\item Form the Markov model with collected transmission probabilities and program deterministic and stochastic simulations
	\item Evaluate coefficient of determination, F-ratio, and t-test statistic for each model
	\item Submit a report evaluating your model. Present (15 minutes) in class.
	\end{itemize}
\end{itemize} 

Students really enjoyed this chance to collect data and be a part of the dataset themselves. Since this personal information was collected anonymously en-masse by an online survey (such as Google Forms without login) and the prescribed intent was to only use it for educational analysis in this class, this project was exempt from IRB (Institutional Review Board) oversight, but I did check with my university's IRB representative. Note, the pertinent exemption is because ``the data can be collected such that individual subjects cannot be identified in any way" \cite{Office_IRB}.

\subsection{Advanced Linear Algebra}
For this theoretical upper division class, I assigned a five-week project analyzing and extending a research article in groups. Having 14-20 students per class, I have used groups of two. Though not an applied project, this is a simple extension of the modeling process to evaluating and improving upon current (undergraduate-level) research articles.

\textbf{Research Article Extension}
\begin{itemize}
\item Resources
	\begin{itemize}
	\item Repository of undergraduate-level mathematics research (e.g. College Mathematics Journal)
	\item Mathematical Typesetting (\href{http://www.lyx.org}{www.lyx.org})
	\item {LyXReport.lyx} and {LyXPresentation.lyx} templates (accessed online)
	\end{itemize}
\item Expectations
	\begin{itemize}
	\item Choose an article on linear algebra theory and replicate its findings
	\item Extend this article by generalizing a theorem or applying the given theory in a new way
	\item Submit a report summarizing the article and detailing your contribution in extending it.
	\item Prepare a 15-minute conference-style presentation
	\end{itemize}
\end{itemize}

Students enjoyed the flexibility to choose their own topic and many of the groups did their best work over the semester on this project. In addition, several of these students have since presented their research projects at regional conferences.

\section{DISCUSSIONS}
In sum, projects are best implemented through incremental communication with students and timely guidance to integrate each project seamlessly into the course. Once the logistical issues are reduced, projects can connect students to course content in exciting new ways and motivate them to learn at a deeper level. 
In talking with a former student who had recently graduated, he noted that my class ``had too much homework, but the projects were great." This particular student was a star pupil and he poured himself into the projects--more hours, I would argue, than he spent on homework during those weeks. Good projects stick in students' minds and infectiously excite future classes. 


\appendix
\addcontentsline{toc}{section}{Appendices}
\section*{APPENDIX}\label{sec:appendix}
\section{PROJECT DETAILS}\label{sec:ProjApp}
This appendix lists the project prompts and scoring guides for each of the projects summarized in section \ref{sec:implementation}.
\subsection{Liberal Arts Mathematics}
\textbf{Monthly Budget}\\
Use the given links to estimate monthly expenses as if you were starting a new job today. The {\tt Monthly Budget Spreadsheet.xlsx} template can be accessed online and uploaded to Google Sheets to become the collaborative spreadsheet for this projects. 

\begin{enumerate}
\item Estimate monthly salary from salary.com and compute taxes withheld (use 11\% Fed + 9\% State (adjust accordingly)) to determine your take-home pay. 
\item Select a home on Zillow to estimate 20\% down due and mortgage payments. Note: it is more realistic that you will rent for a while when you start your job, but it is surprising how quickly you will be interested in buying a house, so it is good to plan ahead for it. How long will it take you to save up enough money for the 20\% down payment. 
\item Check the rate on your student loans and compute the monthly payments needed after you graduate. 
\item Choose how much to put aside for retirement. Use the spreadsheet to calculate how much you will have saved up for retirement in 35 years using an average interest rate of 8\%, and then divide this number by your starting salary's future value (assuming 3\% inflation) to see how many years of living on ``paid vacation" you will have earned. Note: your retirement savings may be augmented by Social Security and retirement accounts your employer sets up for you.
\item Estimate other personal expenses (yellow colored cells) using the example values as a starting point. 
\item Visualize your distribution of Fixed Costs, Financial Goals, and Flexible Spending. How similar is it to the 50-20-30 rule? In what areas could you trim costs to better fit the rule?
\item Alter the amount put into retirement and compare years of retirement and how it impacts your monthly spending. Choose a best value and explain why it is better than the others you checked. Compare with your group members and discuss your differences.
\item Scoring Guide: 100 points total
\begin{itemize}
\item[(20)] Completion of All Components
\item[(20)] Documented Group Collaboration
\item[(20)] Reflection Upon Results/Process
\item[(20)] Overall Presentation of Report
\end{itemize}
\end{enumerate}

\subsection{Discrete Mathematics}
\textbf{Game Analysis}\\
In a group of four, you will choose a simple game for analysis. Past game topics include Planar Graphs, Competing Knight's Tours, Competing N-Queens, Nim, Sim Edge-Coloring, Pipe Layer, Mu Torere, Black Jack, and Battleship. For your chosen simple game, complete the following.
\begin{enumerate}
\item Provide links for further reading and visualizations on this wiki to explain and demonstrate your game to the class.
\item Each person will write up their own strategy (as an algorithm outline) for playing the game in this wiki.
\item Compare strategies through competing multiple times within your group and recording who won the game and in how many moves (or how many possible moves left whichever is easier). Summarize your results on this wiki. Does the first player always win in perfect play or are there conditions on them winning?
\item Lay out a game tree for a three-move end-game scenario and show how the perfect game is selected. Develop an evaluative function which would be a good heuristic for a fixed depth search and demonstrate it.
\item What mathematical concepts are this game based upon?
\item Scoring Guide: 50 points total
\begin{itemize}
\item[(5)] Wiki Construction and Game Visualization
\item[(10)] Outlines of Individual Strategies
\item[(15)] Comparison of Strategies and Tournament Results
\item[(15)] Final Three-Move Game Tree Analysis
\item[(5)] Mathematical Summary of Game
\end{itemize}
\end{enumerate}

\subsection{Differential Equations with Linear Algebra}
Groups of four worked well with the breadth required in these projects for classes between 20 and 43 students. A few semesters when I had 10-14 students, they worked fine with groups of two. The {\tt Common Cold Outbreak Spreadsheet.xlsx} and {Matlab Nutritional Plan.m} files can be accessed online and {\tt Common Cold Outbreak Spreadsheet.xlsx} can be uploaded to Google Sheets as a collaborative spreadsheet.
\subsubsection{Common Cold Outbreak Model}
In groups of four, complete the following activity, develop a model to fit it and write up a wiki report. Due by and presented in class [in week three].
\begin{description}
\item[{Materials}:]  
\begin{itemize}
\item[]
\item Each group gets a printout of a dorm floor plan, 
30 beans, plastic cup, overhead slide piece, and a dry erase marker 
\item Fold along the lines marked and fold in between the two floors to create a `W' shape with the second and third floors down. 
\item Fasten the middle ridge at both ends with paper clips.  
\item Set up a table of time (rounds), number susceptible (not yet infected), number of infected, number of added infected, modeled infected, and modeled infected growth. 
\item Choose a starting infected count between one and 20. Randomly mark that many rooms on a chosen floor with an `X' and set that many beans to the side as infected.
\end{itemize}
\item[{Simulation}:]     
\begin{enumerate} 
\item[] 
\item Run ten simulations where each simulation follows steps 2-5 until there are zero susceptibles or ten rounds have been reached. Use different starting infected counts for each simulation (ranging between one and 20 and distributed randomly on the dorm floor).
\item With one partner walling off the two ends of the chosen floor with their hands, another shakes the susceptible beans from the cup into that floor. Re-shake any beans that fall outside of that floor 
\item Return any beans that do not have a majority of the bean inside an infected room (these are still susceptible and go back in the cup). 
\item For each bean with majority in an infected room, infect the closest available room (mark with an X) and put the beans to the side. Optional: Infect the bathroom only when rooms on both sides are infected and infect hallway when all rooms on at least one side are infected. Otherwise, treat these as ordinary rooms. 
\item Update your table for this round with number of infected and added infected. Number of susceptibles is 30 minus the number infected. 
\item Record all 10 simulations in the class collaborative spreadsheet with 10 tables. Make a combined graph of Round vs. Number Infected for all simulations together.  
\end{enumerate} 


\item[{Analysis}:]
\begin{enumerate} 
\item[]
\item Solve your model and explain the main steps needed in solving it. 
\item For each simulation, estimate the fitting parameter(s) by graphing the general solution with each data set. Show example graphs of your best and worst fitting models over the range of initial infected. 
\item Summarize all 10 values for the fitted parameter(s), their average and range.  
\item Was it a good model? How could this model be improved? 
\item Scoring Guide: 20 points
\begin{itemize}
\item[(5)] Collection/Analysis of Data 
\item[(5)] Analysis/Solution of Model 
\item[(5)] Wiki setup and contribution 
\item[(5)] {Presentation contribution }
\end{itemize}
\end{enumerate}
\end{description}

\subsubsection{Forming a Nutritional Meal Plan}

We will look at RREF, linear independence, and span to select a ``balanced" day of meals from nutrition vectors you will be collecting for your chosen foods. You can either read the nutritional information physically from your food, take a picture of the table at a grocery store, or find the information online. Make sure not to skip values of any of the items listed below in the example. If the \% for a given item is not listed, assume it is zero. You will present your group's balanced day of meals in class [in week 5].
\begin{enumerate}
\item Form a group of four and list your names on the Start page of the Project 2 wiki with your initials inside double ['s to form a link to your page
\item On the bottom of the Start page, collect five food nutrition vectors per person. $A_j$ = [\%Calories (recommended 2000 Calories), \%Total Fat, \%Saturated Fat, \%Cholesterol, \%Sodium, \%Potassium, \%Total Carbohydrates, \%Protein, \%Vit. A, \%Vit. C, \%Calcium, \%Iron, \%Vit. D, \%Vit. E, \%Riboflavin] Food name, Your name, how data was collected

$A_1$ = [70/2000*100, 7, 8, 62, 3, 2, 0, 13, 6, 0, 2, 4, 10, 25, 15]: Egg, Name, data from egg carton (example)
\item Per group, select vectors from the class list and create a matrix, A, which is composed of column vectors $A=[A_1,A_2,...,A_n]$ and b, a column vector of 100's. Using either your calculator, Wolfram Alpha, or Matlab, find rref([A,b]). See the template file {\tt Matlab Nutrition Plan.m} (accessed online). 
\item Create a page that summarizes your balanced day of meals organized into a breakfast, lunch and dinner. For each meal, report the number of servings of each food selected (ignore any negative servings). Using zero for each negative serving, compute your resulting nutrition Ax and discuss which areas your meal plans are excessive/deficient in. Group with most balanced meals (closest to 100\% of each) or most appealing meals get extra credit.
\item Scoring Guide: 20 points total
\begin{itemize}
\item[(5)] Data Collection/Meal Selection
\item[(5)]  Analysis of Data/Results
\item[(10)]  Wiki/Presentation Contribution
\end{itemize}
\end{enumerate}


\subsubsection{Predator-Prey Model}
In groups of four, investigate the interaction between two chosen species. You will give a ten minute summary of your wiki page and models in class [in week 11]. Each group member must contribute to both the wiki and the presentation.

Instructions: 
\begin{enumerate}
\item Choose two species which interact (either predator-prey, competing species, or symbiotic species) 
\item Search for qualitative data (graphs, descriptions of species and their interaction) and quantitative data (average net growth rates, carrying capacities) to set up a 2x2 nonlinear model with parameters. See the following example predator-prey model below.
\begin{eqnarray*}
x' &=& (-ax +by)x\\
y' &=& (-cx + d(1-y))y
\end{eqnarray*}
\item Analyze your model to determine restrictions on the parameters such that all equilibrium points are nonnegative (feasible populations) and that the coexistence equilibrium point is stable or a saddle. 
\item Compute the parameter values using the data by first computing the \% net growth  at two nearby times. Then set up linear equations for a \& b, and c\& d by setting the net growth equal to . If you cannot find quantitative data, you may guess and check parameter values fitting these restrictions and based upon the species' behavior. Type the system into PPLANE applet and plot the solutions. 
\item Summarize your analysis in terms of these populations.
\item Scoring Guide: 30 points total
\begin{itemize}
\item[(10)] Scenario and Model Development with Citations
\item[(10)] Analysis of System 
\item[(10)] Wiki/Presentation Contribution
\end{itemize}
\end{enumerate}

\subsection{Numerical Methods}
The class size for this elective course has fluctuated wildly. I have run this long-term project with individuals for a class of 9, in groups of two for a class of 18, and in groups of four for a class of 31.

\textbf{Final Project}
Choose a topic which can be investigated through a numerical method. Previous topics include electric field of helical solenoid, statistical analysis of voter skepticism, Julia sets (fractals), computing area on the Earth's surface, and heat analysis of air fins.

\begin{description}
\item[Planning:] Sign up for a 30-minute time slot to meet to nail down a detailed proposal for your individual project. Come prepared with your topic and several specific things you would like to investigate about it.

\item[Initialization:] Explain your topic, why you are interested in it, and what kinds of numerical methods would be helpful in investigating your topic (named programs I have given you or ones you have written for this class that would be helpful). Outline three objectives for your project: Easy (takes two hours), Medium (takes one week), Hard (takes more than one week). For each objective, outline a program to complete them in its own Matlab file. Upload to [University's course management system] with this report.

\item[Update:] For each of your Individual Project objectives, continue the list of what you have accomplished, what you still plan to do, and what part(s) you are stuck on/having trouble with. You should be wrapping things up for your final report next Thursday. List any questions you have for me.

\item[Report:]
\begin{enumerate}
\item[]
\item Describe the question you investigated and give some background information on your topic (cite references for all results used).
\item Describe the process used in analyzing the given problem and developing your programs. List all programs used and describe how they are used (including inputs and outputs)
\item Show and describe what each result represents. Then discuss the significance of your results in terms of your investigation and list any avenues that could be interesting to explore further in the future.
\item Scoring Guide: 100 points total
\begin{itemize}
\item[(25)] Typesetting and Organization 
\item[(25)] Understanding of Background Information
\item[(25)] Personal Research Contribution 
\item[(25)] Follows Outlined Objectives
\end{itemize}
\end{enumerate}

\item[Presentation:]
\begin{enumerate} 
\item[]
\item ten-minute presentation of your Final Project
\item Print a one-page flyer for everyone in the class
\item Summarize/demonstrate any code used.
\item Scoring Guide: 50 points total
\begin{itemize}
\item[(10)] Organization and Layout of Presentation 
\item[(15)] Communication of Content
\item[(15)] Evident Contribution and References 
\item[(10)] Engagement of Class \& Ability to Answer Questions
\end{itemize}
\end{enumerate}
\end{description}

\subsection{Mathematical Modeling}
Groups of four worked well with the breadth required in these projects (latest class size was 13). The {\tt Common Cold SIR Spreadsheet.xlsx} and {\tt Cold SIR Regression.m} can be accessed online.
\subsubsection{Hazelnut Estimation}
\begin{enumerate}
\item Estimate the number of hazelnuts that can fit (densely packed) in the jar that I brought to class. The jar is modeled as a large lower cylinder with a smaller upper cylinder on top where the edge retracts to fit the lid. The lower cylinder has diameter 6.6 cm and height 6.5 cm. The upper cylinder has diameter 5.9 cm and height 1.6 cm.
\item Use measurements of your given hazelnut and chosen model representation to represent all hazelnuts in the jar. Estimate the maximum number of hazelnuts that can fit (without crushing) inside the jar. Prizes for best estimate(s).
\item Write a summary explaining your assumptions and approximations used in representing the hazelnut.
\item Scoring Guide: 20 points total
\begin{itemize}
\item[(5)] Description of Model World
\item[(10)] Assumptions and Approximations Stated
\item[(5)] Modeling Process Demonstrated
\end{itemize}
\end{enumerate}

\subsubsection{Common Cold SIR Model}
\begin{enumerate}
\item Over a period of one week, assess the cold symptoms of 40 people (including yourselves). You can either survey people in your living area or send out an online survey. Be sure each person agrees to the use of their data in our analysis beforehand. See example survey form posted. Tabulate the results and compute the number transitioning between states on the collaborative spreadsheet {\tt Common Cold SIR Spreedsheet.xlsx}. 
\item Compute average transmission probabilities between states from your sample data,  form the Markov model, and program the Deterministic model ($X_D$) in Matlab.
\item Use these same probabilities to program a Stochastic Markov model, $X_S$, using the current undergrad population as your total population.
\item Compute $R^2$ values for each of your models to your S, I, R group data, ``Group\_Data". Note, since you just computed averages for your group data (not actual regression), this $R^2 = SS_{Reg}/SS_{Tot}$ actually uses your Markov model populations for $\hat{y}$ not the coefficients. This should be six $R^2$ values, three comparing to the $X_D$ model and three comparing to the $X_S$ model: for example the first $\hat{y}$ is the deterministic model S in $X_D$(1,1:n ) and y is the S data from Group\_Data(1,1:n ), where n is the length of your data vector.
\item  Use the ColdSIR Regression updated file to compute the estimate Markov probabilities for the Class\_Data using multilinear regression. Compare to your group's average values.
\item Compute the F-ratio and t-test statistic of the multilinear regression (b) and evaluate the significance of the trend line and coefficients for the Class\_Data.
\item Write up a report and present your findings (ten minutes) in class. You can do this inside your Matlab file using \%\% as section headers and use the `Publish' command to construct a report with your comments, code, computed results and graphs. See posted ColdSIR Regression file as an example for structuring sections to publish a report.
\item Scoring Guide: 100 points total
\begin{itemize}
\item[(20)] Collection and Presentation of Data
\item[(20)] Statistical Analysis of Data
\item[(25)] Deterministic Model Simulation
\item[(25)] Stochastic Model Simulation
\item[(10)] Discussion
\end{itemize}
\end{enumerate}

\subsection{Advanced Linear Algebra}
\textbf{Research Article Extension}
With a partner, you will find a mathematical idea presented by a mathematical journal article (e.g. College Mathematics Journal) which you can extend in some way. The topic must be relevant to linear algebra theory, not just using matrices in computations. Past project topics include constructing vector spaces, stability of affine transformations, an elementary proof of Dodgson's condensation method for determinants, and generalizing the numerical range of Latin Squares. You will complete a report and oral presentation by [week ten of class].

\begin{description}
\item[Report Guide:]
\begin{enumerate}
\item[] 
\item Demonstrate your knowledge of and work on a topic chosen from a peer-reviewed research article. Detail your contributions with a brief introduction of background information in at least six full pages (two large figures included). 
\item Include necessary definitions, theorems (with at least one proof), and applications or examples to demonstrate your results.
\item Typeset professionally in LyX (or another TeX writer) following the posted {\tt LyXReport.lyx} template.
\item Scoring: 100 pts 
\begin{itemize}
\item[(25)] Typesetting and Organization 
\item[(25)] Understanding of Background Information
\item[(25)] Personal Research Contribution 
\item[(25)] Follows Outlined Research Proposal
\end{itemize}
\end{enumerate}

\item[Presentation Guide:]
\begin{enumerate}
\item[]
\item Oral presentation (10-15 minutes) with slides (Beamer class in LyX is preferred). The {\tt LyXPresentation.lyx} template can be accessed online.
\item Each group member should contribute equally. Start by introducing yourselves. End by thanking the audience. Moderator will ask for questions.
\item Layout should include at least the following slides: Title (title/name/professor/school), 
\item Outline, Intro/Background, Results, Examples/Applications, Discussion/Future Work, Citations/Thank you. 
\item Scoring: 100 pts 
\begin{itemize}
\item[(25)] Organization and Layout of Presentation 
\item[(25)] Communication of Content
\item[(25)] Evident Contribution and References 
\item[(25)] Engagement of Class \& Ability to Answer Questions
\end{itemize}
\end{enumerate}
\end{description}

\section*{BIOGRAPHICAL SKETCH}

{\bf R. Corban Harwood} earned his Ph.D. in Mathematics from Washington State University in 2011 and is currently Assistant Professor of Mathematics at George Fox University in Newberg, OR. In teaching, he loves drawing connections between mathematics and different disciplines such as music, biology, and philosophy. On occasion, Corban consults on modeling projects ranging from testing financial phone apps to designing optimal battery chemistries for hybrid-electric vehicles. His research interests include semilinear partial differential equations, oscillatory analysis of numerical methods, and modeling reaction-diffusion phenomena. On the side, he enjoys cycling in the Willamette Valley and hiking to waterfalls with his family.


\end{document}